\documentstyle[12pt]{article}
\newtheorem{thm}{Theorem}
\newtheorem{lem}{Lemma}
\newtheorem{rem}{Remark}

\newtheorem{ex}{Example}
\newtheorem{coro}{Corollary}

\begin{document}

\normalsize

\begin{center}
{\LARGE\bf
Polynomials with Real Zeros and\\
P\'olya Frequency Sequences}
\\[10pt]
{Yi\ Wang$^{\rm a}$
\footnote
{Partially supported by NSF of Liaoning Province of China
Grant No. 2001102084}
\quad and\quad
Yeong-Nan\ Yeh$^{\rm b}$
\footnote
{Partially supported by NSC
92-2115-M-001-016}}
\\[10pt]
$^{\rm a}$ {\footnotesize Department of Applied Mathematics, Dalian
University of Technology, Dalian 116024, China}
\\[5pt]
$^{\rm b}$ {\footnotesize Institute of Mathematics, Academia Sinica,
Taipei 11529, Taiwan}
\\[5pt]
\end{center}
\vspace{20pt}
\centerline{\bf Abstract} \vspace{7pt}

Let $f(x)$ and $g(x)$ be two real polynomials
whose leading coefficients have the same sign.
Suppose that $f(x)$ and $g(x)$ have only real zeros
and that $g$ interlaces $f$ or $g$ alternates left of $f$.
We show that if $ad\ge bc$ then the polynomial
$$(bx+a)f(x)+(dx+c)g(x)$$
has only real zeros.
Applications are related to certain results of F.Brenti
(Mem. Amer. Math. Soc. 413 (1989))
and transformations of P\'olya frequency sequences.
More specifically,
suppose that $A(n,k)$ are nonnegative numbers
which satisfy the recurrence
$$
A(n,k)=(rn+sk+t)A(n-1,k-1)+(an+bk+c)A(n-1,k)
$$
for $n\ge 1$ and $0\le k\le n$,
where $A(n,k)=0$ unless $0\le k\le n$.
We show that if $rb\ge as$ and $(r+s+t)b\ge (a+c)s$,
then for each $n\ge 0$, $A(n,0),A(n,1),\ldots,A(n,n)$
is a P\'olya frequency sequence.
This gives a unified proof of the PF property of many well-known sequences
including the binomial coefficients,
the Stirling numbers of two kinds and the Eulerian numbers.
\\[7pt]
\\[7pt]

\section{Introduction}
\hspace*{\parindent}
Let $a_0,a_1,a_2,\ldots$ be a sequence of nonnegative real numbers.
It is {\em unimodal} if
$a_0\le a_1\le\cdots\le a_{k-1}\le a_k\ge a_{k+1}\ge\cdots$
for some $k$.
It is {\em log-concave} (LC)
if $a_{i-1}a_{i+1}\le a_i^2$ for all $i>0$.
It is said to {\em have no internal zeros}
if there are no three indices $i<j<k$ such that
$a_i,a_k\not=0$ and $a_j=0$.
Clearly, a log-concave sequence
with no internal zeros is unimodal.
Unimodal and log-concave sequences occur naturally
in combinatorics, algebra, analysis, geometry, computer science,
probability and statistics.
We refer the reader to the survey papers
by Stanley\cite{Sta89} and Brenti\cite{Bre94}
for various results on unimodality and log-concavity.

A classical approach for attacking unimodality and log-concavity problems of finite sequences
is to use the following result
originally due to Newton\cite[p.104]{HLP}.
\\[7pt]
{\bf Newton's Inequality}\quad
{\em
Given a finite sequence $a_0,a_1,\ldots,a_n$ of nonnegative numbers.
Suppose that its generating function
$\sum\limits_{i=0}^na_ix^i$
has only real zeros.
Then
$$
a_i^2\ge a_{i-1}a_{i+1}
\left({i+1\over i}\right)\left({n-i+1\over n-i}\right),\qquad i=1,2,\ldots,n-1,
$$
and in particular,
the sequence is log-concave(with no internal zeros).
}
\vspace{7pt}

It is natural to look at those sequences
whose associated polynomial has only real zeros.
A characterization for such sequences
comes from the theory of total positivity.
Let $A=(a_{ij})_{i,j\ge 0}$ be an infinite matrix.
We say that $A$ is {\em totally positive} (or {\em TP}, for short)
if all minors of $A$ have nonnegative determinants.
An infinite sequence $a_0,a_1,a_2,\ldots$ of nonnegative numbers is called
a {\em P{\'o}lya frequency sequence} (or a {\em PF} sequence)
if the matrix $(a_{i-j})_{i,j\ge 0}$ is a TP matrix
(where $a_k=0$ if $k<0$).
A finite sequence $a_0,a_1,\ldots,a_n$ is {\em PF}
if the infinite sequence $a_0,a_1,\ldots,a_n,0,0,\ldots$ is PF.
By definition,
a PF sequence is necessarily log-concave (with no internal zeros).
A deeper result is the following theorem
which provides the basic link between
finite PF sequences and polynomials having only real zeros.
For more information about TP matrices and PF sequences,
see \cite{Kar68}.
\\[7pt]
{\bf Aissen-Schoenberg-Whitney Theorem}\cite{ASW52}
\quad
{\em
A finite sequence $a_0, \ldots, a_n$ of nonnegative numbers is PF
if and only if its generating function
$\sum\limits_{i=0}^na_ix^i$
has only real zeros.
}
\vspace{1pt}

Brenti was the first one who applied total positivity techniques
to study systematically unimodality and log-concavity problems\cite{Bre89}.
He noted that PF sequences have much better behavior
than unimodal and log-concave sequences.
It may often be more convenient to prove that a sequence is PF
even if we are actually interested only in the unimodality or log-concavity.
On the other hand,
many unimodal and log-concave sequences arising in combinatorics
turn out to be PF sequences.
So it is natural to pay more attention to PF sequences,
as well as polynomials having only real zeros.
Indeed, there are a number of open problems in combinatorics
concerning whether certain polynomials have only real zeros
(see \cite{Bre94, Sta89, Sta00} for instance).
The present paper is devoted to studying certain transformations of polynomials
(resp. sequences) that preserve the reality of zeros (resp. PF property).

Let ${\bf RZ}$ be the set of real polynomials having only real zeros
and ${\bf PF}$ the set of polynomials in ${\bf RZ}$
whose coefficients are nonnegative.
In other words,
${\bf PF}$ is the set of polynomials whose coefficients form a PF sequence.
It is clear that all zeros of each polynomial in ${\bf PF}$ are non-positive.
For convenience let ${\bf RZ}$ contain all real constants
and ${\bf PF}$ contain all nonnegative numbers.

Suppose that $f(x)\in{\bf RZ}$ and $g(x)\in{\bf RZ}$.
Let $r_n\le\cdots\le r_1$ and $s_m\le\cdots\le s_1$
be the zeros of $f$ and $g$ respectively.
Following \cite{Wag92},
we say that $g$ {\em interlaces} $f$ if $m=n-1$ and
$$r_n\le s_{n-1}\le r_{n-1}\le \cdots\le r_2\le s_1\le r_1,$$
and that $g$ {\em alternates left of} $f$ if $m=n$ and
$$s_n\le r_n\le s_{n-1}\le r_{n-1}\le \cdots\le r_2\le s_1\le r_1.$$
By $g\leadsto f$ we denote
``either $g$ interlaces $f$ or $g$ alternates left of $f$".
For notational convenience,
let $a\leadsto  bx+c$ for all real constants $a,b,c$.
Clearly, $g\leadsto  f$ yields $ag\leadsto bf$ for any $ab\ne 0$.

The main result of this paper is the following.
\begin{thm}\label{int}
Let $f(x)$ and $g(x)$ be two real polynomials
whose leading coefficients have the same sign
and let $F(x)=(bx+a)f(x)+(dx+c)g(x)$.
Suppose that $f,g\in{\bf RZ}$ and $g\leadsto f$.
Then, if $ad\ge bc$, $F(x)\in{\bf RZ}$.
\end{thm}

This paper is organized as follows.
The next section is devoted to the proof of Theorem \ref{int}.
In \S 3 we present some applications of Theorem \ref{int}
related to certain results of Brenti\cite{Bre89} and transformations of PF sequences,
the latter of which induces a unify proof of the PF property of many well-known sequences,
including the binomial coefficients,
the Stirling numbers of two kinds, and the Eulerian numbers.

\section{Proof of Theorem \ref{int}}
\hspace*{\parindent}
Let $sgn$ denote the sign function defined by
$$
{\rm sgn} (t)
=\cases{+1 & for $t>0$,\cr \ \ 0 & for $t=0,$\cr -1 & for $t<0$.}
$$
Let $f(x)$ be a real function.
If $f(x)>0$ (resp. $f(x)<0$) for sufficiently large $x$,
then we denote ${\rm sgn} f(+\infty)=+1$(resp. $-1$).
The meaning of ${\rm sgn} f(-\infty)$ is similar.

Before showing Theorem \ref{int},
we provide three lemmas to deal with certain special cases of the theorem.
The first one is a fundamental and well-known result
(see Section 3 of \cite{Wag92} for instance).
\begin{lem}\label{blem}
Suppose that $f, g\in{\bf RZ}$ and $g\leadsto f$.
Then $f+g\in{\bf RZ}$.
Furthermore,
if the leading coefficients of $f$ and $g$ have the same sign,
then $g\leadsto f+g\leadsto f$.
\end{lem}

Let $P(x)$ be a real polynomial of degree $n$.
Define its {\em reciprocal} polynomial by
$$P^*(x)=x^nP(1/x).$$
The following facts are elementary but very useful in the sequel:

(i)\quad
If $P(0)\ne 0$ then $\deg P^*=\deg P$ and $\left(P^*\right)^*=P$.

(ii)\quad
$P\in{\bf RZ}$ if and only if $P^*\in{\bf RZ}$.

(iii)\quad
Suppose that all zeros of $f$ and $g$ are negative and $g\leadsto f$.
If $g$ interlaces $f$, then $g^*$ interlaces $f^*$.
If $g$ alternates left of $f$, then $f^*$ alternates left of $g^*$.

In the following two lemmas
we assume that $f$ and $g$ are two monic polynomials with only simple negative zeros
and that $g$ interlaces $f$.
More precisely, let
$f(x)=\prod\limits_{i=1}^n(x-r_i)$ and
$g(x)=\prod\limits_{i=1}^{n-1}(x-s_i)$
where
$$r_n<s_{n-1}<r_n<\cdots<s_2<r_2<s_1<r_1<0$$
($g(x)=1$ provided $n=1$).
\begin{lem}\label{ad>0}
Suppose that $ad>0$. Then

{\em(i)}\quad
$af+dxg\in{\bf RZ}$ and $f,g\leadsto  af+dxg$;

{\em(ii)}\quad $af+(dx+c)g\in{\bf RZ}$ for any $c$;

{\em(iii)}\quad $(bx+a)f+dxg\in{\bf RZ}$ for any $b$.
\end{lem}
{\bf Proof}\quad Let $F=af+dxg$. Then $F^*=af^*+dg^*$.

(i)\quad
By the assumption,
$g$ interlaces $f$, so $g^*$ interlaces $f^*$.
From Lemma \ref{blem} it follows that $F^*\in{\bf RZ}$,
$g^*$ interlaces $F^*$ and $F^*$ alternates left of $f^*$.
Thus $F\in{\bf RZ}$, $g$ interlaces $F$ and $f$ alternates left of $F$.

(ii)\quad
By (i), $F\in{\bf RZ}$ and $g\leadsto  F$.
Thus for any $c$, $F+cg\in{\bf RZ}$ by Lemma \ref{blem}, i.e.,
$af+(dx+c)g\in{\bf RZ}$.

(iii)\quad
By (i), $f\leadsto  F$.
Note that all zeros of the polynomial $F$ are negative
since its coefficients have the same sign.
Hence $F\leadsto  xf$.
Thus for any $b$, $bxf+F\in{\bf RZ}$ by Lemma \ref{blem},
i.e., $(bx+a)f+dxg\in{\bf RZ}$. \qquad$\Box$
\begin{rem}
{\em
The condition $ad>0$ in Lemma \ref{ad>0} is not necessary.
Actually, $g\leadsto  f$ implies $f\leadsto  xg$
since all zeros of $f$ are negative,
and thus for any $a$ and $d$,
$af+dxg\in{\bf RZ}$ by Lemma \ref{blem}.
We can also show that $ac>0$ implies
$af+(dx+c)g\in{\bf RZ}$ for any $d$,
which is in a sense ``dual'' to Lemma \ref{ad>0}(ii)
since $[af+(dx+c)g]^*=af^*+(cx+d)g^*$.
Similarly, as the dual version of Lemma \ref{ad>0}(iii),
we have $(bx+a)f+dxg\in{\bf RZ}$ when $bd>0$.
}
\end{rem}
\begin{lem}\label{bc<0}
Suppose that $bc<0$. Then

{\em(i)}\quad
$bxf+cg\in{\bf RZ}$ and $f,xg\leadsto  bxf+cg$;

{\em(ii)}\quad $bxf+(dx+c)g\in{\bf RZ}$;

{\em(iii)}\quad $(bx+a)f+cg\in{\bf RZ}$
and $f\leadsto  (bx+a)f+cg$ for any $a$.
\end{lem}
{\bf Proof}\quad
Let $F=bxf+cg$.
Without loss of generality, let $b>0$ and $c<0$.
We have
$$
{\rm sgn} F(r_j)=
{\rm sgn} \left(c\prod\limits_{i=1}^{n-1}(r_j-s_i)\right)=(-1)^j,
\quad j=1,2,\ldots,n
$$
and
$$
{\rm sgn} F(s_j)=
{\rm sgn}
\left(bs_j\prod\limits_{i=1}^n(s_j-r_i)\right)
=(-1)^{j+1}, \quad j=1,2,\ldots,n-1.
$$
Also,
${\rm sgn} F(0)=-1, {\rm sgn} F(-\infty)=(-1)^{n+1}$ and ${\rm sgn} F(+\infty)=1$.

By the intermediate-value theorem,
$F$ has $n+1$ real zeros
$t_1,\ldots,t_{n+1}$ satisfying
$$t_{n+1}<r_n<s_{n-1}<t_n<r_{n-1}\cdots<s_2<t_3<r_2<s_1<t_2<r_1<0<t_1.$$
Thus $F\in{\bf RZ}$ and $f,xg\leadsto  F$.
This proves (i).

Write $bxf+(dx+c)g=F+d(xg)$ and$(bx+a)f+cg=F+af$.
Then both (ii) and (iii) follow from (i) and Lemma \ref{blem}.
\qquad$\Box$
\\[7pt]
{\bf Proof of Theorem \ref{int}}\quad
Without loss of generality,
we may assume that $f$ and $g$ are monic
and have no zeros in common
(which implies that they have only simple zeros).
We may also assume that all zeros of $f$ and $g$ are negative.
Since if we define $f_1(x)=f(x+u), g_1(x)=g(x+u)$ and $F_1(x)=F(x+u)$
where $u$ is a real number larger than all zeros of $f$ and $g$,
then $f_1(x), g_1(x)$ have only negative zeros and $g_1\leadsto f_1$.
Moreover,
$$F_1(x)=(b_1x+a_1)f_1(x)+(d_1x+c_1)g_1(x)$$
where $a_1=bu+a,b_1=b,c_1=du+c,d_1=d$.
Clearly, $ad\ge bc$ is equivalent to $a_1d_1\ge b_1c_1$ and
$F(x)\in{\bf RZ}$ is equivalent to $F_1(x)\in{\bf RZ}$.
Thus we may consider $F_1$ instead of $F$.

Suppose first that $g$ interlaces $f$.
If $f(x)=x-r$ and $g(x)=1$,
then
$$F(x)=bx^2+(a-br+d)x-(ar-c).$$
The discriminant of $F(x)$ is
$$(a-br+d)^2+4b(ar-c)=(a+br-d)^2+4(ad-bc)\ge 0.$$
Thus $F(x)\in{\bf RZ}$. Now let
$f(x)=\prod\limits_{i=1}^n(x-r_i)$
and $g(x)=\prod\limits_{i=1}^{n-1}(x-s_i)$ where $n>1$ and
$$r_n<s_{n-1}<r_n<\cdots<s_2<r_2<s_1<r_1<0.$$

If $abcd=0$,
then the statement follows from the previous three lemmas.
If $ad=bc$,
then the statement follows from Lemma \ref{blem}
since $F$ has the factor $bx+a$.
So, let $abcd\ne 0, ad>bc$ and $a>0$.
We distinguish four cases.

{\bf Case 1.}\quad $c>0$.
Then $g\leadsto  af+cg$ by Lemma \ref{blem}.
Note that
$$aF=(bx+a)(af+cg)+(ad-bc)xg.$$
Hence $aF\in{\bf RZ}$ by Lemma \ref{ad>0}(iii),
and so $F\in{\bf RZ}$.

{\bf Case 2.}\quad $b,c<0$.
Then $d>0$ and $-a/b>-c/d>0$.
Thus $(dx+c)g$ interlaces $(bx+a)f$,
and so $F=(bx+a)f+(dx+c)g\in{\bf RZ}$ by Lemma \ref{blem}.

{\bf Case 3.}\quad $b,d>0$.
Then $g\leadsto  bf+dg$ by Lemma \ref{blem}.
Note that
$$bF=(bx+a)(bf+dg)+(bc-ad)g.$$
Hence $bF\in{\bf RZ}$ by Lemma \ref{bc<0}(iii),
and so $F\in{\bf RZ}$.

{\bf Case 4.}\quad $b>0,c,d<0$.

Assume that $r_k<-c/d<r_{k-1}$ for some $k$.
Then $dr_j+c<0$ for $1\le j\le k-1$ and $dr_j+c>0$ for $k\le j\le n$.
We have
$$
{\rm sgn} F(r_j) ={\rm
sgn}\left((dr_j+c)\prod_{i=1}^{n-1}(r_j-s_i)\right) =\cases{(-1)^j
& if $1\le j\le k-1$,\cr (-1)^{j+1} & if $k\le j\le n$,}
$$
and
$$
{\rm sgn} F\left(-{c\over d}\right)
= {\rm sgn}\left({ad-bc\over d}\prod\limits_{i=1}^n
\left(-{c\over d}-r_i\right)\right)
=(-1)^{k}.
$$
Also, ${\rm sgn} F(+\infty)=1$.
Thus $F(x)$ has $n+2$ changes of sign.
This implies that $F(x)$ has $n+1$ real zeros.
So $F\in{\bf RZ}$.

Similarly,
we may prove $F\in{\bf RZ}$
provided $r_1<-c/d$ or $r_n>-c/d$.

Assume now that $r_k=-c/d$ for some $k$.
We consider only the case $1<k<n$
since the proof for the case $k=1$ or $k=n$ is similar.
We have
$$
{\rm sgn} F(r_j) ={\rm
sgn}\left((dr_j+c)\prod_{i=1}^{n-1}(r_j-s_i)\right) =\cases{(-1)^j
& if $1\le j\le k-1$,\cr 0 & if $j=k$,\cr (-1)^{j+1} & if $k+1\le
j\le n$,}
$$
and ${\rm sgn} F(+\infty)=1$.
Also,
$$
{\rm sgn} F(s_k)
={\rm sgn}\left((bs_k+a)\prod\limits_{i=1}^n(s_k-r_i)\right)=(-1)^{k+1}
$$
since $bs_k+a<br_k+a={ad-bc\over d}<0$.

Thus $F(x)$ has $n-k$ zeros in the interval $(r_n,s_k)$
and $k-1$ ones in the interval $(r_{k-1},+\infty)$.
It remains to show that the interval $(s_k,r_{k-1})$ contains two zeros of $F(x)$.
Actually, since $F(s_k)$ and $F(r_{k-1})$ have the same sign,
the interval $(s_k,r_{k-1})$ contains an even number of zeros of $F(x)$
(see \cite[Part V, Prob.8]{PS} for instance).
In this interval, $F(x)$ has one zero $-c/d$ and thus at least two zeros,
as desired.

Next suppose that $g$ alternates left of $f$. Let
$f(x)=\prod\limits_{i=1}^n(x-r_i)$ and
$g(x)=\prod\limits_{i=1}^{n}(x-s_i)$ where
$s_n<r_n<\cdots<s_1<r_1$. Define
\begin{eqnarray*}
f_1(x)&=&\prod\limits_{i=1}^n(x+r_1-s_i)=g(x+r_1),\\
g_1(x)&=&\prod\limits_{i=2}^{n}(x+r_1-r_i)=f(x+r_1)/x
\end{eqnarray*}
(set $g_1(x)=1$ if $n=1$),
and $F_1(x)=F(x+r_1)$.
Then $f_1,g_1$ have only negative zeros and $g_1$ interlaces $f_1$.
Moreover,
$$F_1=[dx+(dr+c)]f_1+x[bx+(br+a)]g_1.$$
So
$$F_1^*=[(dr+c)x+d]f_1^*+[(br+a)x+b]g_1^*.$$
Clearly, $\deg f_1^*=n, \deg g_1^*=n-1$ and $g_1^*$ interlaces $f_1^*$.
Also,
$$d(br+a)=bdr+ad\ge bdr+bc=(dr+c)b.$$
By the result of the first part,
we have $F_1^*\in{\bf RZ}$.
Hence $F_1\in{\bf RZ}$, and so $F\in{\bf RZ}$.

Thus the proof of the theorem is complete.\qquad$\Box$

\section{\large\bf Applications of Theorem \ref{int}}
\hspace*{\parindent}
In this section we give some applications of Theorem \ref{int}.
In \cite{Bre89},
Brenti investigated linear transformations that preserve the PF property.
Our results are closely related to those of Brenti.
The first one is an immediate consequence of Theorem \ref{int}.
\begin{coro}\label{abcd}
Suppose that $f(x), g(x)\in{\bf PF}$ and that $g$ interlaces $f$.
Let
$$F(x)=(ax+b)f(x)+x(cx+d)g(x).$$
If $ad\ge bc$, then $F(x)\in{\bf RZ}$.
\end{coro}
{\bf Proof}\quad
Note that $g$ interlaces $f$ implies $f$ alternates left of $xg$
since all zeros of $f$ are non-positive.
Thus the statement follows from Theorem \ref{int}.
\qquad$\Box$

Corollary \ref{abcd} generalizes Theorem 2.4.4 of Brenti\cite{Bre89}
which states $F(x)\in{\bf RZ}$ when $ac<0$ and $ad>\max\{bc, 0\}$.
Let we consider a special case of Corollary \ref{abcd}.
Recall that if $f(x)\in{\bf PF}$ then $f'(x)\in{\bf PF}$ and $f'$ interlaces $f$.
By Corollary \ref{abcd}
we have that
$(ax+b)f(x)+x(cx+d)f'(x)\in{\bf RZ}$
provided $ad\ge bc$.
This result can be restated in terms of sequences
instead of polynomials and their derivatives.
\begin{coro}\label{derseq}
Let $x_0,x_1,\ldots,x_{n-1}$ be a PF sequence and
$$y_k=[a+c(k-1)]x_{k-1}+(b+dk)x_k,\quad k=0,1,\ldots,n,$$
where $x_{-1}=x_n=0$.
If $ad\ge bc$ and all $y_k$ are nonnegative,
then the sequence $y_0,y_1,\ldots,y_n$ is PF.
\end{coro}

Corollary \ref{derseq} actually gives a class of linear transformations
that preserve the PF property of sequences.
It may also be used to reproduce certain results of Brenti, e.g.,
\cite[Theorem 2.4.2]{Bre89}.

We next use Corollary \ref{derseq} to study the PF property of triangular arrays.
Let $\{x(n,k)\}_{n\ge k\ge 0}$ be a triangular array of nonnegative numbers
satisfying a two-term recursion
$$
x(n,k)=b(n,k)x(n-1,k-1)+a(n,k)x(n-1,k)
$$
for $n\ge 1, 0\le k\le n$,
where $x(n,k)=0$ unless $0\le k\le n$,
and the coefficients $a(n,k),b(n,k)$ are nonnegative.
For convenience we take $x(0,0)=1$.
Such a triangular array is said to be {\em unimodal}(resp. {\em LC}, {\em PF})
if for each $n\ge 0$,
the sequence $x(n,0),x(n,1),\ldots,x(n,n)$ has the corresponding property.
Canfield\cite{Can84} considered the unimodality property of triangular arrays.
Kurtz\cite{Kur72} and Sagan\cite{Sag88}
provided certain sufficient conditions on coefficients such that
triangular arrays are LC respectively.
What conditions will insure that triangular arrays are PF?
This question is, in general, very difficult to answer
(see Brenti\cite[Theorem 4.3]{Bre95} for partial solution).
Since our interest in this matter stems from combinatorial motivations,
we consider only triangular arrays $\{A(n,k)\}$
which satisfy a recurrence of ``bilinear" form
\begin{eqnarray}\label{bilin}
A(n,k)=(rn+sk+t)A(n-1,k-1)+(an+bk+c)A(n-1,k)
\end{eqnarray}
for $n\ge 1, 0\le k\le n$,
where $A(n,k)=0$ unless $0\le k\le n$, and $r,s,t,a,b,c$ are real numbers.
Many important triangular arrays arising in combinatorics satisfy such a recurrence.
It is easy to see that such triangular arrays are necessarily LC\cite{Kur72}.
By means of Corollary \ref{derseq},
we may obtain the following.
\begin{coro}\label{pftri}
Let $\{A(n,k)\}$ be a triangular array defined by (\ref{bilin}).
Suppose that $rb\ge as$ and $(r+s+t)b\ge (a+c)s$.
Then the triangular array $\{A(n,k)\}$ is PF.
\end{coro}
{\bf Proof}\quad
We need to show that for each $n\ge 0$,
the sequence $A(n,0),\ldots,A(n,n)$ is PF.
We proceed by induction on $n$.
The result is clearly true for $n=0,1$.
So suppose that $n>1$ and $A(n-1,0),\ldots,A(n-1,n-1)$ is PF.
By the assumption we have
\begin{eqnarray*}
(rn+s+t)b-(an+c)s&=&(rb-as)n+(s+t)b-cs\\
&\ge& (rb-as)+(s+t)b-cs\\
&=&(r+s+t)b-(a+c)s\\
&\ge& 0.
\end{eqnarray*}
Thus $A(n,0),\ldots,A(n,n)$ is PF by Corollary \ref{derseq},
and the proof is therefore complete.
\qquad$\Box$

In what follows we list some examples of triangular arrays
whose PF property has been proved by various techniques in the literature.
All of these triangular arrays are easily seen to
satisfy the assumption of Corollary \ref{pftri}.
So the PF property is an immediate consequence of Corollary \ref{pftri}.
\begin{ex}
{\em
The binomial coefficients $n\choose k$,
the (signless) Stirling numbers of the first kind $c(n,k)$,
the Stirling numbers of the second kind $S(n,k)$
and the classical Eulerian numbers $A(n,k)$
(the numbers of permutations of $1,2,\ldots, n$
having $k-1$ descents)
satisfy the recurrence
\begin{eqnarray*}
{n\choose k}&=&{n-1\choose k-1}+{n-1\choose k},\\
c(n,k)&=&c(n-1,k-1)+(n-1)c(n-1,k),\\
S(n,k)&=&S(n-1,k-1)+kS(n-1,k),\\
A(n,k)&=&(n-k+1)A(n-1,k-1)+k A(n-1,k)
\end{eqnarray*}
respectively (see \cite{Sta86} for instance).
It is well known that generating functions of ${n\choose k}$ and $c(n,k)$
are $(x+1)^n$ and $x(x+1)\cdots(x+n-1)$, respectively.
So the binomial coefficients and
the (signless) Stirling numbers of the first kind are clearly PF.
Harper\cite{Har67} showed the PF property of the Stirling numbers of the second kind
(see \cite{Dob68,KL78,Kur72,Lie68} for log-concavity).
For the PF property of the Eulerian numbers, see \cite[p.292]{Com74}.
}
\end{ex}
\begin{ex}
{\em
The associated Lah numbers defined by
$$L_m(n,k)=\left(n!/k!\right)\sum_{i=1}^k(-1)^{k-i}{k\choose i}{n+mi-1\choose n}$$
satisfy the recurrence
$$L_m(n,k)=mL_m(n-1,k-1)+(mk+n-1)L_m(n-1,k).$$
Ahuja and Enneking\cite{AE79} showed that
$\{L_m(n,k)\}_{0\le k\le n}$ is PF for all $m\ge 1$.
}
\end{ex}
\begin{ex}
{\em
The associated Stirling numbers of two kinds,
introduced by Jordan and Ward
(see \cite{Car71} for instance),
satisfy the recurrence
$$c^*(n,k)=(2n-k-1)\left(c^*(n-1,k-1)+c^*(n-1,k)\right)$$
and
$$S^*(n,k)=(n-k)S^*(n-1,k-1)+(2n-k-1)S^*(n-1,k)$$
respectively.
Kurtz\cite{Kur72} pointed out both of them are LC,
and Ahuja\cite{Ahu73} showed that they are PF.
}
\end{ex}
\begin{ex}
{\em
The holiday numbers $\psi(n,k)$ and $\phi(n,k)$ of the first kind and the second kind,
introduced by L. A. Sz\'ekly in \cite{Sze85}, satisfy the recurrence
$$\psi(n,k)=\psi(n-1,k-1)+(2n+k-1)\psi(n-1,k)$$
and
$$\phi(n,k)=\phi(n-1,k-1)+(2n+k)\phi(n-1,k)$$
respectively.
Sz\'ekly\cite{Sze87} showed that both of them are PF.
}
\end{ex}
\begin{ex}
{\em
It is a long-standing conjecture that
the Whitney numbers of the second kind of any finite geometric lattice are unimodal
or even log-concave sequence (see \cite{Wel76}).
Dowling\cite{Dow73} constructed a class of geometric lattices
over a finite group of order $m$
and showed that the Whitney numbers $W_m(n,k)$ satisfy the recurrence
$$W_m(n,k)=W_m(n-1,k-1)+(1+mk)W_m(n-1,k).$$
Stonesifer\cite{Sto75} showed the log-concavity of $W_m(n,k)$.
Benoumhani proved that $W_m(n,k)$ and $k!W_m(n,k)$ are PF
in \cite{Ben99} and \cite{Ben97} respectively.
It is easy to see that $k!W_m(n,k)$ actually satisfies the recurrence
$$A(n,k)=kA(n-1,k-1)+(1+mk)A(n-1,k).$$
}\end{ex}

It is worth mentioning that certain examples to
which our results are not applicable can be made to be so
by means of appropriate transformations.
Let us examine a result of R. Simion
which concerns the unimodality of the numbers $\Theta ({\bf n},k)$
of compositions of a multiset ${\bf n}=(n_1,n_2,\ldots)$ into exactly $k$ parts.
It is known that the following recurrence holds
(see \cite[p.96]{Rio78}):
\begin{eqnarray}\label{sim1}
(n_j+1)\Theta ({\bf n}+e_j,k)
=k\Theta ({\bf n},k-1)+(n_j+k)\Theta ({\bf n},k),
\end{eqnarray}
where ${\bf n}+e_j$ denotes the multiset obtained from ${\bf n}$
by adjoining one (additional) copy of the $j$th type element.
Let $f_{{\bf n}}(x)=\sum\limits_{k\ge 0}\Theta ({\bf n},k)x^k$
be the corresponding generating function.
Then by (\ref{sim1}),
\begin{eqnarray}\label{sim2}
(n_j+1)f_{{\bf n}+e_j}(x)
=(x+n_j)f_{{\bf n}}(x)+x(x+1)f'_{{\bf n}}(x).
\end{eqnarray}
Simion showed that the multiplicity $m$ of $-1$ as a zero of
$f_{{\bf n}}(x)$ is $\max\limits_{i}\{n_i-1\}$
(see \cite{Sim84} for details).
Based on this and the recurrence (\ref{sim2}),
she showed that $f_{{\bf n}}(x)\in{\bf RZ}$ implies
$f_{{\bf n}+e_j}(x)\in{\bf RZ}$.
(Thus $f_{{\bf n}}(x)\in{\bf RZ}$ for any multiset ${\bf n}$
and $\Theta ({\bf n},k)$ is therefore unimodal in $k$.)
Note that this result can not be derived directly from Theorem \ref{int}
since $n_j\ge 1$. However,
if we write $f_{{\bf n}}(x)=(x+1)^mg_{{\bf n}}(x)$,
then the problem can be reduced to show that
$g_{{\bf n}}(x)\in{\bf RZ}$ implies
$[(m+1)x+n_j]g_{{\bf n}}(x)+x(x+1)g'_{{\bf n}}(x)\in{\bf RZ}$.
This follows immediately from Corollary \ref{abcd} since $m+1\ge n_j$.

\section{Concluding Remarks}
\hspace*{\parindent}
There are various methods for showing that
polynomials have only real zeros.
One basic method for showing that a sequence $P_0(x),P_1(x),P_2(x),\ldots$ of polynomials
has real zeros is to show by induction
that the polynomials have interlaced simple (real) zeros,
that is, they form a {\em Sturm sequence}\cite{Sta89}.
The key part of this method is in the inductive step.
The sequences of polynomials occurred in combinatorics often satisfy certain recurrence relations.
Theorem \ref{int} provides a useful tool for solving this kind of problems.
For example, let $\{p_n(x)\}$ be a sequence of orthogonal polynomials associated with certain distribution.
Then the following three-term recursion holds:
\begin{eqnarray}
p_n(x)=(a_nx+b_n)p_{n-1}(x)-c_np_{n-2}(x),\quad n=1,2,\ldots,
\end{eqnarray}
with $p_{-1}(x)=0$,
where $a_n, b_n$ and $c_n$ are real constants and $a_n,c_n>0$
(see \cite[Theorem 3.2.1]{Sze39} for instance).
Well-known examples of orthogonal polynomials
include the Jacobi polynomials,
the Hermite polynomials and the Laguerre polynomials.
By induction and Lemma \ref{bc<0}(iii),
it follows immediately that
$p_n(x)$ has only real and simple zeros
and $p_n$ interlaces $p_{n+1}$ for $n\ge 1$,
which is a well-known result.

Many special cases of Theorem \ref{int} have occurred in the literature.
We refer the reader to, for example,
\cite{GW96,HL72,Stahl97,Wag91,Wag92}.
\vspace{20pt}

\noindent \centerline{\bf Acknowledgments}

This research was completed during the first author's stay in the
Institute of Mathematics, Academia Sinica, Taipei. The first
author would like to thank the Institute for its support.

The authors would like to thank the anonymous referee for his/her
careful reading and valuable suggestions that led to an improved
version of this manuscript.

\end{document}